\documentclass{article}
\usepackage{amsmath, amsfonts, amsthm, graphicx, amssymb}

\newtheorem{theorem}{Theorem}
\newtheorem{lemma}{Lemma}
\newtheorem{corollary}{Corollary}
\newtheorem{claim}{Claim}
\newtheorem{defn}{Definition}
\newtheorem{conjecture}{Conjecture}

\newcommand{\opsp}{\mbox{\ensuremath{S^{1} \setminus P(G)}}}
\newcommand{\adda}[1]{\ensuremath{#1_{\alpha}}}

\newcommand{\orb}[1]{\ensuremath{\mathcal{O}(#1)}}

\newcommand{\orbh}[1]{\ensuremath{\mathcal{O}_{H^{1}}(#1)}}
\newcommand{\rg}[1]{\ensuremath{\mathcal{R}(#1)}}
\newcommand{\pfbox}{\begin{flushright}$\Box$\end{flushright}\vspace{1pt}}

\begin{document}
\title{Orbit structure for groups of homeomorphisms of $S^{1}$}
\author{N. C. Esty \thanks{ncesty@gmail.com}}
\maketitle

\begin{abstract}
In this paper we will refine Sacksteder's theorem for groups of orientation-preserving homeomorphisms of the circle in the case that there exists a finite orbit set.  We will give a categorization of the topological possibilities for the orbits of points of the circle, along with examples.
\end{abstract}

\section{Introduction}
\label{sec:intro}

In this paper we consider groups of orientation-preserving
homeomorphisms of the circle, acting in the natural way, and give a
complete categorization of the topological possibilities for orbits of
points of the circle under the action of the group. This is an
extension of Sacksteder's theorem (\cite{ghys}, p. 20), which gives
three possible orbit types, but does not describe the remaining
possibilities.  In particular, if there is a finite orbit set, it does
not describe the structure of the remaining orbits.

In Section~\ref{sec:prelim} we will review the Poincar\'{e}
classification of orbits for circle homeomorphisms, and recall
Sacksteder's theorem for groups of homeomorphisms.  The intent of this
paper is to refine Sacksteder's theorem in the case that the action of
the group has a finite orbit.

In Section~\ref{sec:groups} we will prove the following:

\begin{lemma} 
Let $G$ be a group of orientation-preserving homeomorphisms of the
interval $\bar{I} = [0,1]$ and suppose $G$ satisfies one of the following
conditions:
\begin{enumerate}
\item $G$ is finitely generated and there are no fixed points of the
group in $(0,1)$.
\item There exists at least one group element $g$ of $G$ whose fixed
points do not accumulate at either 0 or 1.
\end{enumerate}
Then every nonempty $G$-invariant set that is closed with respect to
$I = (0,1)$ contains a nonempty, $G$-invariant set closed in $(0,1)$
and minimal under those properties.
\end{lemma}

We will assume throughout the rest of the paper that the group $G$ of
orientation-preserving homeomorphisms of the circle has a finite
number of finite orbit points, which we denote $P(G)$.  We will split
the circle at points of $P(G)$ and look at the complementary
subintervals, which we call $\adda{C}$.  If the finite orbits points
in $P(G)$ are all fixed under the group, then the orbit of any point
will be entirely contained in one $\adda{C}$; if not, it will be
contained in the union of some number of the $\adda{C}$.  We call the
union of the intervals which the orbit of $x$ intersects nontrivially
$\rg{x}$, for \emph{range} of $x$.  We assume the group $G$, restricted to
one $\adda{C}$, satisfies one of the conditions of the lemma.  We then
prove the following:

\begin{theorem}
For all $x \in S^{1} \setminus P(G)$, there are four possibilities for
the orbit of $x$:
\begin{enumerate}
\item $\orb{x}$ is dense in $\rg{x}$.
\item $\orb{x}$ accumulates exactly at points of $P(G)$.
\item $\orb{x}$ is contained in a $G$-invariant Cantor set in
$\overline{\rg{x}}$.  The orbit of $x$ is dense in the Cantor set, and
the Cantor set is contained in the closure of the orbit of each point
of $\rg{x}$.
\item $\overline{\orb{x}}$ contains a proper subset which is
$G$-invariant, closed with respect to $S^{1} \setminus P(G)$, and
which does not intersect $\orb{x}$; i.e., $\orb{x}$ accumulates at the
closure of another orbit.
\end{enumerate}
\end{theorem}

In Section~\ref{sec:examples} we will give a collection of examples to
demonstrate the four cases of the theorem.  In particular, we will
give two different examples that have a Cantor set, and we will look
at several examples designed to explore the different topological
possibilites of a point in Case 4.  This will lead us to
define the \emph{level} of an orbit (see \ref{sec:leveln}), so that we
can describe the nesting of orbits that may occur in this case.

In Section~\ref{sec:semigroups} we will consider semigroups rather
than groups, and give an example which demonstrates one way in which
the orbit type under a semigroup can differ from those for groups.

Section~\ref{sec:case4} is devoted to a closer examination of Case 4,
the case in the theorem that admits the most complexity of orbits.  We
will prove a collection of lemmas whose purpose is to reduce the
number of considerations necessary to understand the range of
possibilities, and we will define a few of the special orbit types
that may arise.  

In Section~\ref{sec:analytic} we will prove results for points in the
fourth case when the group homeomorphisms are analytic.

Finally, in Section~\ref{sec:qsandconjs} we will conclude with some
conjectures for future work.

\section{Preliminaries}
\label{sec:prelim}

The classical theory of dynamical systems studies the orbit structure
of a homeomorphism or a flow on a manifold, in particular, the
topological properties of $\{f^{n}(x) : n \in \mathbb{Z}\}$ or
$\{\phi(x, t): t \in \mathbb{R}\}$.  This corresponds to
studying actions of the group $\mathbb{Z}$ or $\mathbb{R}$.  One can
ask under what circumstances a map will have periodic points, dense
orbits, etc.  One can also study the actions of more general groups.
In this paper we will look at the orbit structure of a group of
homeomorphisms acting in a natural way on the circle.

Poincar\'{e} classified the possible behavior of orbits for a given
homeomorphism of the circle by looking at the rotation number of the
homeomorphism, an invariant that gives the ``average'' rotation under
the map.  For a more thourough discussion of these concepts, see
\cite{katok}, Chapter 11.  He showed that points of $S^{1}$ can have
six different kinds of orbits under homeomorphisms: three types each
for rational and irrational rotation numbers.  The possibilities are:
a periodic orbit, a homoclinic orbit approaching a given periodic
orbit, a heteroclinic orbit approaching two different periodic orbits,
a dense orbit, an orbit dense in a Cantor set, and an orbit homoclinic
to a Cantor set.

When studying the action of a group of homeomorphisms, we discuss the
\emph{group orbit} of a point $x$, denoted $\orb{x} = \{\phi(g)(x) : g
\in G\}$.  We often refer to this as the \emph{orbit} of the point
without making reference to the group.  This is not to be confused
with the orbit of a point under a single map (generally a ``smaller''
set).  Sacksteder's theorem describes the main possibilities for the
dynamics of an arbitrary group of homeomorphisms.  See \cite{ghys} for
a discussion and proof.

\begin{theorem}[Sacksteder] Let $G$ be any subgroup of
$\mathbf{Homeo}_{+}(S^{1})$, the orientation-preserving homeomorphisms
of the circle.  There are three mutually exclusive possibilities:
\begin{enumerate}
\item  There is a finite orbit.
\item All orbits are dense.
\item There is a compact $G$-invariant subset $K \subset S^{1}$
which is infinite and different from $S^{1}$ and such that the orbits
of points in $K$ are dense in $K$.  This set $K$ is unique, contained
in the closure of any orbit, and is homeomorphic to a Cantor set.
\end{enumerate}
\end{theorem}

The three options are directly related to the types of invariant
minimal sets which can exist for homeomorphisms of the circle.  It is
known that the circle supports three different types of invariant
minimal sets under the action of a homeomorphism: finite sets, Cantor
sets, and the circle itself.  As the orbit of a point under a
homeomorphism is automatically invariant under that map, it is easy to
see how the results of the classification theorems reflect those
options.  Unsurprisingly, the proof of this theorem relies heavily on
the ability to find a minimal set which is compact and invariant under
the action of the group.  We are interested in extending the result
when there is a known finite orbit, in order to understand the
behavior of all the nonfinite orbits in this case.  We will also need
to find a minimal set; however, it will take some extra work.

\section{Groups of Homeomorphisms}
\label{sec:groups}

Unless indicated otherwise, throughout the rest of this paper we shall
assume the following:

\begin{enumerate}

\item $G < \mathbf{Homeo}_{+}(S^{1})$ acts on the circle in the
natural way, i.e., for all $g \in G$, $\phi(g)(x) = g(x)$.

\item The action of $G$ on the circle has a (nonzero)
finite number of finite orbits. Denote by $P(G)$ be the set of points of
finite orbits.

\item $G$ is either finitely generated, or has at least one group
element with isolated fixed points.  

\end{enumerate}

Notationally, we will use $I = (0,1)$ to be the open unit
interval, and $\bar{I} = [0,1]$ to be the closed.  

\begin{defn} If a set $K$ is
invariant under the action of the group $G$ and closed with respect to
a proper subset $S \subset S^{1}$, we will say that $K$ is a
\emph{$G_{S}$-set}.  If $K$ contains no proper subset that is also
$G$-invariant and closed in $S$, we will say that $K$ is
\emph{$G_{S}$-minimal}.
\end{defn}

We use the following lemma repeatedly:

\begin{lemma}

Let $G$ be a group of orientation-preserving homeomorphisms of
$\bar{I}$.  Suppose either of the following
conditions are satisfied:
\begin{itemize}
\item[\textbf{C1}] There exists at least one group element whose fixed points
do not accumulate at 0 and do not accumulate at 1.
\item[\textbf{C2}] G is finitely generated and there are no fixed points of
the group in $(0, 1)$.
\end{itemize}
Then each nonempty $G_{I}$-set $M$ contains a nonempty set $K$ which is
$G_{I}$-minimal.
\label{lemma:zorn}
\end{lemma}

\emph{Proof:} The existence of a compact minimal set is well known
when the space is compact.  We will find a closed interval $J$ in
$(0,1)$, depending only upon the group $G$, and show that any
$G_{I}$-set $M$ must intersect $J$.  We then consider $\mathcal{M}$,
the collection of all nonempty $G_{I}$-sets.  We order $\mathcal{M}$
by downward inclusion and apply Zorn's lemma to produce a maximal
chain $\{M_{\alpha}\}$.  The intersection $\bigcap M_{\alpha}$ belongs
to $\mathcal{M}$.

First, suppose $G$ satisfies C1.  Let $g$ be the specified group
element and choose $J$ to be a nonempty closed interval in $(0,1)$
containing all the fixed points of $g$ in its interior, large enough
that $g(J) \cap J \neq \emptyset$; we write $J = [a, b]$.

Let $M \in \mathcal{M}$.  Assume $x \in M$, but $x \notin J$; the
situation is symmetric, so without loss of generality suppose that $0
< x < a$.  The point $x$ is not fixed by $g$, so (by switching to
$g^{-1}$ if neccesary) we may assume that if $0 < x < a$ then $x <
g(x)$.  We know $a < g(a)$ and $g(a) < b$.  Therefore positive
iterates of $g$ can not ``jump over'' $J$.  For all $n$, $g^{n}(x)$
lies in $M$.  If the sequence $(g^{n}(x))_{n \in \mathbb{N}}$
accumulates at a point $y$, then $y$ is a fixed point of $g$ and must
be interior to $J$.  For large enough $n$, $g^{n}(x) > a$.  Therefore
$M \cap J \neq \emptyset$ and we have shown that Condition 1 is
sufficient to imply the conclusion of the lemma.

Now suppose $G$ satisfies C2.  Let $\{g_{1}, \dots, g_{k}\}$ be the
generators of $G$.  Choose $J$ to be a nonempty closed interval in
$(0,1)$ large enough that for all $i$, \mbox{$g_{i}(J) \cap J \neq
\emptyset$}, and write $J = [a,b]$.  Suppose $M$ is a $G_{I}$-set, and
take $x \in M, x \notin J$; without loss of generality, suppose $x <
a$.

Since there are no fixed points for the group inside $(0,1)$, for
every $x$ there is at least one generator $g$ with $g(x) \neq x$.  As
we can always switch to $g^{-1}$, it is fine to assume that between
two consecutive fixed points of a generator $g_{i}$, $g_{i}(x) > x$.
We would like to find a sequence of generators that will move $x$
steadily to the right, and so into $J$; the concern is that the only
such sequences of generators, when applied in order to form the
sequence $(x_{n})_{n \in \mathbb{N}}$, $x_{n} = g_{n}(x_{n-1})$, will
create a sequence of points that converges to some $y \notin J$.
Suppose this is the case.  If this is problematic, then $y < a$.  If,
in some $\epsilon$-neighborhood of $y$, there are no fixed points of
one particular generator $g$, then once our sequence of points enters
that $\epsilon$-neighborhood we may apply $g$ repeatedly and the
resulting sequence of points will eventually move past $y$.  Since we
assumed no sequence of generators would do this, all neighborhoods
must contain fixed points of all the generators.  This means for each
$i$, $0 \leq i \leq k$, there is a sequence of fixed points of $g_{i}$
converging to $y$: $p^{i}_{n} \rightarrow y$.  But then, by
continuity, $y$ is fixed for $g_{i}$.  As this is true for each $i$,
$y$ is a fixed point of the group, and we have a contradiction.
Therefore some sequence of generators must cause $x_{n}$ to be in $J$
for large enough $n$, and so $M \cap J \neq \emptyset$.  We have shown
that Condition 2 is sufficient.

\pfbox

\begin{corollary}  
If $G$ is a group of analytic homeomorphisms of $\bar{I}$, then any
$G_{I}$-set $M$ contains a nonempty set $K$ that is $G_{I}$-minimal.
\end{corollary}

The proof is omitted.

Let \{$\adda{C}$\} be the finite collection of connected components of
$S := \opsp$.  Let $C_{\alpha(x)}$ be that component which contains
$x$.  We define the following set, the \emph{range} of $x$:
\[ \rg{x} = \bigcup_{\adda{C} \cap \orb{x} \neq \emptyset} \adda{C}  = \mathcal{O}(C_{\alpha(x)})\] 
Note that if the finite orbit points are all fixed under the group,
$\rg{x} = C_{\alpha(x)}$.  The following is clear.

\begin{lemma} 
If $x$, $y$ are members of the same component $\adda{C}$, then $\rg{x} = \rg{y}$.
\end{lemma}

We can now state the main theorem:

\begin{theorem}
\label{thm:main}
Let $G < \mathbf{Homeo}_{+}(S^{1})$ be a finitely generated group with
a finite number of finite orbit points $P(G)$.  For every $x \in \opsp
$, there are four possibilities for the orbit of $x$,
$\mathcal{O}(x)$:
\begin{enumerate}
\item $\orb{x}$ is dense in $\rg{x}$.
\item $\orb{x}$ accumulates only at points of $P(G)$.
\item $\orb{x}$ is contained in a $G$-invariant Cantor set in
$\overline{\rg{x}}$.  The orbit of $x$ is dense in the Cantor set, and
the Cantor set is contained in the closure of the orbit of each point
of $\rg{x}$.
\item $\overline{\orb{x}}$ contains a proper subset which is
$G$-invariant, closed with respect to \opsp, and which does not
intersect $\orb{x}$; i.e., $\orb{x}$ accumulates along the closure of
another orbit.
\end{enumerate}
\end{theorem}

\emph{Remark:} If the group is not finitely generated, it is enough
that for each $\adda{C}$ the group contains an element whose fixed
points do not accumulate at the endpoints of $\adda{C}$.  We merely
require that in each $\adda{C}$ the group satisfies either condition
of Lemma~\ref{lemma:zorn}.

\emph{Proof:}

Consider the collection of subsets of $S := \opsp$ which are
$G_{S}$-sets.  For any $x \in S$, $\overline{\orb{x}} \cap S$ is a
$G_{S}$-set.  We take a maximal chain of such sets starting with
$\overline{\orb{x}} \cap S$; as $G$ is finitely generated, and we are
looking at the group elements restricted to a particular $C_{\alpha}$
(topologically equivalent to $I$), Condition 1 of the lemma is
satisfied, and there exists a $G_{S}$-minimal element $K$.  (Again, if the
group is not finitely generated, but satisifes Condition 2 for
$\adda{C}$, we also have a minimal set.)  Note that different maximal
chains may exist, as the ordering is partial, and so there may be many
such minimal sets corresponding to different maximal chains in
$\overline{\orb{x}}$.  The proof breaks into two sections: the first
considers points that are members of some such $G_{S}$-minimal set,
and the second considers points that are not.

First suppose $x \in K$.  Because $K$ is closed with respect to $S$
and invariant under $G$, $K$ must contain $\overline{\orb{x}} \cap
S$. By minimality, $K = \overline{\orb{x}} \cap S$, and so all points
of $K$ have orbits dense in $K$.

Consider the subsets of $K$, $\partial K = K \setminus \mbox{Int}(K)$,
the boundary of $K$, and $K' \cap S$,\ the set of accumulation
points of $K$ inside $S$.  Both are easily shown to be $G_{S}$-sets.
Because $K$ is $G_{S}$-minimal, this allows for three mutually
exclusive possibilities for $\orb{x}$:

\begin{enumerate} 

\item $\partial K = \emptyset$, in which case $K$ is equal to its
interior, and is therefore open with respect to $S$.  However, it
was closed with respect to $S$ by construction, and so it must
be the union of some number of the connected components $\adda{C}$.
Since $\orb{x}$ is dense in $K$, $K = \rg{x}$.

\item $K' \cap S = \emptyset$. The orbit of $x$ has an infinite number
of points and $S^{1}$ is compact, so the set of accumulation points
of the orbit can not be empty.  If it is empty in $S$, $K'$ must be
contained in $P(G)$.  So the orbit of $x$ accumulates only at points
of $P(G)$.

\item If $\partial K$ and $K' \cap S$ are not empty, then by
minimality, $\partial K = K = K' \cap S$.  This means $K$ has no
interior and all points of $K$ are accumulation points.
$\overline{K}$ is therefore perfect, nonempty, closed, and has no
interior; it is a Cantor set.  All points of $\overline{K}$ that are
not in $P(G)$ are in $K$, and have dense orbits in $K$.  It is clear
that $\overline{K} \subset \overline{\rg{x}}$.

\end{enumerate}

The fourth case of the theorem arises from the only remaining
possibility: that the point $x \in S^{1}$ is not contained in any of
the $G_{S}$-minimal sets $K$ obtained by intersecting a maximal chain
of $G_{S}$-sets.  Consider such a maximal chain of $G_{S}$-sets in
$\overline{\orb{x}}$.  Since the corresponding $G_{S}$-minimal set
does not contain $x$, neither can it contain any points of the orbit
of $x$.  In order for this to occur,

\begin{enumerate}
\item[4.] $\overline{\orb{x}}$ must have a proper subset, $K$, that is
also a $G_{S}$-set, and this subset must be contained entirely in
$\overline{\orb{x}} \setminus \orb{x}$.
\end{enumerate}

Therfore $\orb{x}$ must fall into one of the four cases in the
theorem.  We now show the Cantor set $\overline{K}$ of Case 3 is
contained in the closure of the orbit of each point of $\rg{x}$:

Let $y$ be a point of $\rg{x}$.  If $y \in K$, then we already know
the orbit of $y$ is dense in $K$, so assume $y \in \rg{x} \setminus
K$.  The complement of $K$ in $\rg{x}$ is a countable collection of
open intervals, so $y$ is in the interior of some $I$.  At least one
endpoint of $I$ must be in $K$: call that endpoint $a$.

Let $p \in K$.  Since $a \in K$, there exist group elements $g_{n}$
such that $(g_{n}(a))$ is a sequence of distinct points that converges
to $p$.  Because $(g_{n}(a))$ is a distinct sequence, $g_{n}(I)$ must
be an infinite collection of intervals; the size of the intervals must
go to zero, and so $d(g_{n}(a), g_{n}(y)) \rightarrow 0$.  Since
$g_{n}(a) \rightarrow p$, $g_{n}(y) \rightarrow p$.  Therefore the
orbit of $y$ is dense in $K$.

\pfbox

As each $\adda{C}(x)$ is part of $\rg{x}$, the orbit of $x$ in one
\adda{C} will be topologically equivalent to the orbit of $x$ in another.
Therefore we will assume that all points of $P(G)$ are fixed.  This
allows us to switch to homeomorphisms of the closed unit interval
$\bar{I} = [0,1]$, where 0 and 1 will correspond to two consecutive
members of $P(G)$. 

\section{Examples}
\label{sec:examples}

Let us suppose that we have a point $y$ that falls into Case 4 of the
theorem: \orb{y} accumulates along the closure of another orbit, say
that of the point $z$.  We want to describe \orb{y}.

If $z$ is in Case 1, namely $\orb{z}$ is dense in $I$, then clearly
the orbit of $y$ is also dense in $I$.  If $z$ is in Case 4, then we
could pick a point from Case 1, 2 or 3 on whose orbit it accumulated,
and \orb{y} would accumulate there as well.  We assume that $z$ falls
into either Case 2 (which we will call \emph{integer type}, due to the
fact that there is a homeomorphism taking $I$ to $\mathbb{R}$ and
taking \orb{z} to the integers) or into Case 3, the Cantor set case.

If $z$ is either integer type or Cantor set type, then the complement
of the closure of the orbit of $z$ consists of a countable number of
open intervals $I_{n}$, such that $y \in I_{0}$.  If $\orb{z}$ has
integer type, the intervals will be adjacent and we use the natural
notation $I_{i} = [z_{i}, z_{i+1}]$ (where $z = z_{0}$).  Often in
this case there is a particular group element $f$ with $f^{i}(z_{0}) =
z_{i}$, as we shall see in the examples.  If $\orb{z}$ is a Cantor
set, we will use the notation $I_{i} = [z_{i}^{1}, z_{i}^{2}]$ (where
$z = z_{0}^{1}$), because the $I_{i}$ are not adjacent and do not
follow a natural left-to-right ordering.

\subsection{Case 1: Dense Orbits}

Let $f = x^{1/3}$ and $g = x^2$.  Let $G$ be the group generated by
$f$ and $g$.  Clearly $\orb{x} = \{x^{2^{k}/3^{j}} | k, j \in
\mathbb{Z}\}$.  Because $\ln(2)/\ln(3)$ is irrational, this is a dense
set.  This example is equivalent to the translation of $\mathbb{R}$ by
rationally independent numbers.

\subsection{Case 2: $\orb{x}$ Accumulates Only at Points of $P(G)$}
 
If the group $G$ is generated by a single homeomorphism $g$, then for
any non-fixed point $x$, the orbit will consist of a countable
collection of points, the iterates of $x$ under the map $g$.  This set
will accumulate exactly at two consecutive fixed points of $g$, and so
$x$ is in Case 2.  Later we will refer to this as a \emph{level 1
integer type} orbit; see Section~\ref{sec:leveln} for the definition
of \emph{level}.

\subsection{Case 3: Cantor Set Type, Example 1}

For Case 3 we have two examples: one that uses an infinitely countable
number of generators, and one that uses only two.  We first look at
the former, as the concept is simpler.  In this case, our group will
satisfy the second condition of Lemma~\ref{lemma:zorn}, as it is not
finitely generated. We use the following lemma, a slight variation on
the fact that Cantor sets in the real line are ambiently homeomorphic.
See \cite{pugh} for a proof.

\begin{lemma}
\label{lemma:cantor1}
Given two Cantor sets in $\mathbb{R}$, $C_{1}$ and $C_{2}$, and given
for each a point in the Cantor set which is a left endpoint of one of
the open intervals making up the complement of the set, $l_{1}$ and
$l_{2}$, then there exists a homeomorphism from $\mathbb{R} \rightarrow
\mathbb{R}$ that takes $C_{1}$ to $C_{2}$ and $l_{1}$ to $l_{2}$.
\end{lemma}

Let $C$ be the standard middle thirds Cantor set in $[0,1]$.  Let
$\{x_{n}\}$ be the countable collection of left endpoints of $C$.

We ``split'' the Cantor set at a given left endpoint into two Cantor
sets, one to each side of the interval which the endpoint borders.  We
do this for $x_{0}$ and $x_{1}$.  Using the lemma, let $f_{0}$ be the
homeomorphism which maps the Cantor set to the left of $x_{0}$ to the
Cantor set to the left of $x_{1}$, and the Cantor set to the right of
$x_{0}$ to the Cantor set to the right of $x_{1}$, and takes the
interval $I_{0}$ to $I_{1}$ in a linear manner.  This map preserves
$C$ and takes $x_{0}$ to $x_{1}$.

Continuing in the same manner, let $f_{1}$ be the homeomorphism that
takes $x_{1}$ to $x_{2}$ and preserves $C$, and generally, $f_{n}$
that homeomorphism which takes $x_{n}$ to $x_{n+1}$ and preserves $C$.
The result is a countable number of generators, each of which preserve
$C$.  Clearly the orbit of $x_{0}$ contains all the $x_{n}$ and is
therefore dense in $C$, and so $\orb{x}$ is Case 3, Cantor Set Type.

\subsection{Case 3: Cantor Set Type, Example 2}

In order to construct a Case 3 example using only two generators, it
is important to have a clear picture of the Cantor set.  To ease the
discussion, we start by setting up some notation.

\subsubsection{Notation}

Let $\bar{I} = [0,1]$ and let $C$ be the standard middle thirds Cantor set.
Let $g$ be the following piecewise linear homeomorphism of $\bar{I}$
preserving $C$:

\[ g(x) = \left\{ \begin{array}{ll}
  3x & \mbox{if $0 \leq x \leq \frac{2}{9}$} \\
  x + 4/9 & \mbox{if $\frac{2}{9} < x < \frac{1}{3}$} \\
  \frac{1}{3}(x-1) + 1 & \mbox{if $\frac{1}{3} < x \leq 1$}
  \end{array}
\right. \]

\begin{figure}
\resizebox{\textwidth}{!}{\includegraphics{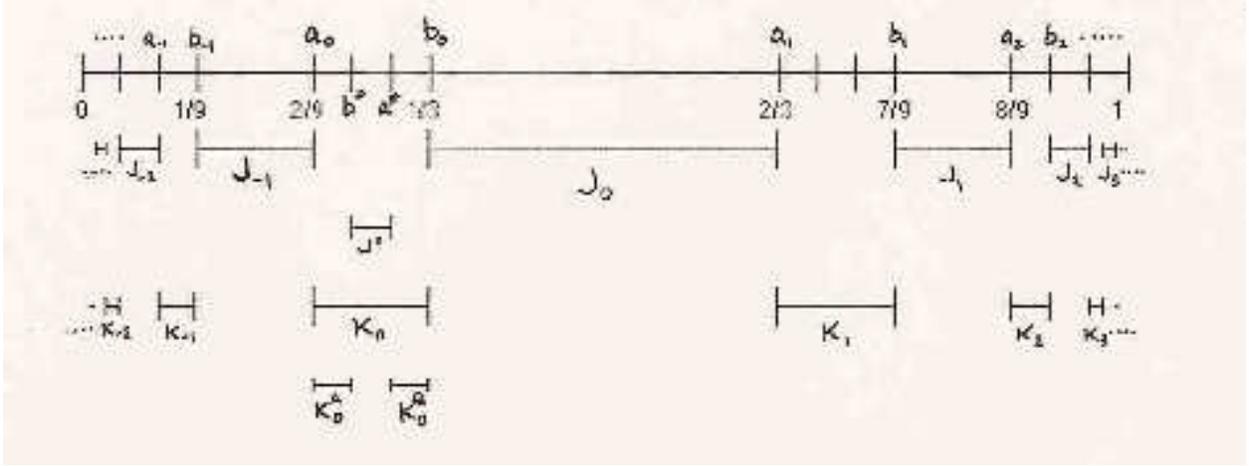}}
\caption{Notation for the Cantor set case.}
\label{fig:cantornotation}
\end{figure}

The notation is pictured in Figure~\ref{fig:cantornotation}.

Let $K_{0} = [2/9, 1/3]$, $J_{0} = (1/3, 2/3)$, $a_{0} = 2/9$, and
$b_{0} = 1/3$.  Considered independently, $K_{0}$ contains its own
middle thirds Cantor set, $C_{0} = C \cap K_{0}$.  We use the
following notation: For all $n \in \mathbb{Z}$, $K_{n} = g^{n}(K_{0}),
J_{n} = g^{n}(J_{0}), a_{n} = g^{n}(a_{0})$, and $b_{n} =
g^{n}(b_{0})$.  Each $K_{n}$ is a closed interval containing a middle
thirds Cantor set $C_{n}$, which is the image of $C_{0}$ under
$g^{n}$.  This will allow us to easily indicate which section of the
Cantor set contains a given point.

We could instead have defined $a_{n}$ as the intersection of the
closures of $J_{n-1}$ and $K_{n}$: accordingly, we will sometimes
refer to $a_{n}$ as the ``left edge'' of $C_{n}$ or the ``right
endpoint'' of $J_{n-1}$.  Similarly, we could think of $b_{n}$ as the
intersection of the closures of $K_{n}$ and $J_{n}$: in that respect,
$b_{n}$ is the ``right edge'' of $C_{n}$ or the ``left endpoint'' of
$J_{n}$.

\emph{Remark:} The open unit interval $I$ is the disjoint union of the
$K_{n}$ and $J_{n}$, and the standard Cantor set $C$ is the disjoint
union of the points $\{0\}$, $\{1\}$, and each $C_{n}$.  We describe
points $x \in C$ as belonging to one of three categories: if $x = 0$
or $x = 1$, we say that $x$ is a \emph{terminal edge}.  If $x = a_{n}$
or $x = b_{n}$ for some $n$, we say $x$ is an \emph{edge}.  Otherwise
we say that $x$ is \emph{interior to $C$}.  This is a slight abuse of
terminology, as technically, $C$ has no interior, but it is helpful is
describing the relationship between $x$ and other points of the Cantor
set.

Given a point $x \in C$, we associate to it its ternary expansion,
which we call its \emph{position}.  Let $P$ be the set of all
positions of points in the Cantor set, and $P^{L}$ the set of
positions of those $x \in C$ that are left endpoints of one of the
connected components of the complement of $C$ (namely, the $b_{n}$).
Although $C$ (and therefore the set of all positions) is uncountable,
$P^{L}$ is countable.  We put an ordering $<$ on $P$ (and $P^{L}$) in
the natural way: $p_{1} < p_{2}$ if the point given by the ternary
expansion $p_{1}$ is less than the point given by the expansion
$p_{2}$.

Since each $C_{n}$ is also a middle thirds Cantor set, for any
$x \in C$ other than $0$ and $1$, we can describe $x$ by specifying
which subinterval contains it and its position in the Cantor set of
that subinterval: $x = (K_{n}, p)$.  We will use this notation
interchangeably with $x$ for points in the Cantor set.

\emph{Remark:} $g^{k}(K_{n}) = K_{n+k}$, and similarly for $J_{n}$,
$a_{n}$ and $b_{n}$.  Because $g$ is linear on any given $K_{n}$, it
preserves position: i.e., if $x = (K_{n}, p)$, then $g^{k}(x) =
(K_{n+k}, p)$.

Finally, we define one further subdivision of $K_{0}$: Let $K_{0}^{A}
= [2/9, 7/27]$, $J^{*} = (7/27, 8/27)$, $K_{0}^{B} = [8/27, 1/3]$, and
let $b^{*} = 7/27, a^{*} = 8/27$.  Like the other $K_{n}$, $K_{0}^{A}$
and $K_{0}^{B}$ both contain middle third Cantor sets, $C_{0}^{A}$ and
$C_{0}^{B}$.  Using our (interval, position) notation, this means
\[ \begin{array}{lclcl}
  b^{*} & = & (K_{0}^{A},\, 1) & = & (K_{0},\, .1) \\
   a^{*} & = & (K_{0}^{B},\, 0) & = & (K_{0},\, .2) 
   \end{array}  \]

However, this last subdivision is only necessary in $K_{0}$.  

\subsubsection{Construction of f}

We consider the set $D = \{(p_{1}, p_{2}, p^{*}_{1}, p^{*}_{2}):
p_{i}, p^{*}_{i} \in P^{L}, p_{1} < p_{2}, p^{*}_{1} < p^{*}_{2}\}$,
of quadruples of positions. Because $P^{L}$ is countable, so is $D$.
Let $\pi:D \rightarrow \mathbb{Z}_{+}$ be a one-to-one and onto map
from the quadruples to the positive integers.  We will use the following
variation on Lemma~\ref{lemma:cantor1}.

\begin{lemma}
Given two middle thirds Cantor sets $C \subset \bar{I}$ and $C^{*}
\subset \bar{I^{*}}$, and given any $(p_{1}, p_{2}, p_{1}^{*},
p_{2}^{*}) \in D$, there exists a homeomorphism $f: \bar{I}
\rightarrow \bar{I^{*}}$ such that $f$ takes $C$ onto $C^{*}$ and
(using our alternate notation for points of $C$), $f(\bar{I}, p_{1}) =
(\bar{I^{*}}, p_{1}^{*})$ and $f(\bar{I}, p_{2}) = (\bar{I^{*}},
p^{*}_{2})$.  In other words, $f$ preserves the Cantor set and sends one
pair of left endpoints to the other pair.
\end{lemma}

Let $n \geq 1$.  Since $\pi$ is an onto function, $n = \pi((p_{1},
p_{2}, p^{*}_{1}, p^{*}_{2}))$ for some quadruple in $D$.  We define $f$
on $K_{n}$ as the homeomorphism from $K_{n}$ to $K_{n+1}$ given by the
lemma, meaning $f$ takes $C_{n}$ to $C_{n+1}$ and takes the pair of
left endpoints $(p_{1}, p_{2})$ to $(p^{*}_{1}, p^{*}_{2})$.
Let $f$ take $J_{n}$ to $J_{n+1}$ linearly.

For $x \in \bar{I}$ with $x \leq a_{-1}$, let $f$ equal $g$.  In other
words, for $n < -1$, $f(K_{n}) = K_{n+1}$ and $f$ preserves position
in each $C_{n}$.

\begin{figure}
\resizebox{\textwidth}{!}{\includegraphics{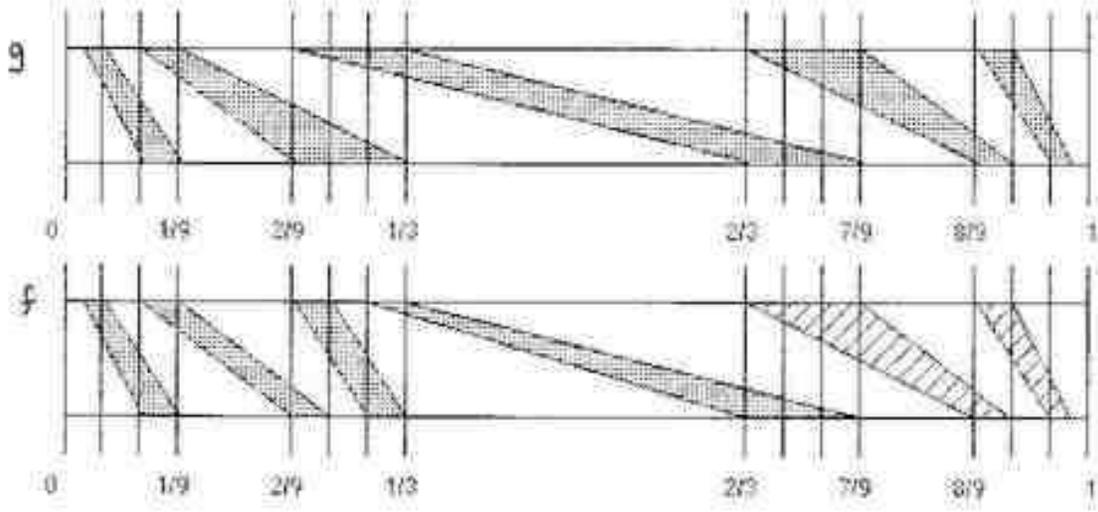}}
\caption{The difference between the images of $K_{i}$ under $f$ and $g$.
Dotted areas indicate that the map preserves position,
cross-hatched areas indicate that positions are not preserved.}
\label{fig:fvsg}
\end{figure}

Let $f(K_{-1}) = K_{0}^{A}$, taking $C_{-1}$ to $C_{0}^{A}$ and
preserving position; let $f(K_{0}^{A}) = K_{0}^{B}$, taking
$C_{0}^{A}$ to $C_{0}^{B}$ and preserving position; finally, let
$f(K_{0}^{B}) = K_{1}$, taking $C_{0}^{B}$ to $C_{1}$ and preserving
position.  As in the rest of $I$, $f$ maps the $J$ intervals linearly.
This completes the definition of $f$ on all of $\bar{I}$: See
Figure~\ref{fig:fvsg}.

\emph{Note:} If $n \geq 1$ and $k \geq 1$, or if $n \leq -1$ and $k
\leq 0$, then $f^{k}(K_{n}) = K_{n+k}$.

\subsubsection{Cantor Set Example 2}

Let $G$ be the group of homeomorphisms of $\bar{I}$ generated by $f$
and $g$.  For $x \in I$ we have, as always, $\orb{x} = \{ h(x) \, : \,
h \in G\}$.

\begin{theorem} 
If $x \in C$ and $x \neq 0$, $x \neq 1$, then $\overline{\orb{x}} = C$.
\end{theorem}

This implies that under the group $G$, points of $C \cap I$ fall into
Case 3: orbits which are dense in a Cantor set.  Recall the three
categories for points of $C$: terminal edge, edge, and interior.  We
will use these terms to describe different cases that the proof needs
to address.  

\emph{Proof:}

Let $x \in C$, and $x \neq 0, 1$.  Let $y \in C$ and let $\epsilon >
0$.  The point $y$ may be an edge, a terminal edge, or interior to
$C$, whereas the point $x$ can not be a terminal edge.  For all cases,
we will demonstrate the existence of a group element $h \in G$ such
that
\begin{equation}
|h(x) - y| < \epsilon   
\label{eq:claim}
\end{equation}

Suppose first that $y$ is a terminal edge.  If $y = 1$, then there
exists some $n$ sufficiently large that $1 - g^{n}(x) < \epsilon$.  If
$y = 0$, then there exists some $n$ sufficiently large that $g^{-m}(x)
< \epsilon$.  So assume that $y$ is either an edge point, or interior
to $C$.  We write $y = (K_{m}, p_{m})$ and $x = (K_{n}, p_{n})$.

If $p_{n} = 1$, namely, if $x$ is a right edge point, we define the
map $h_{1} = gf^{-1}g^{-n}$.  If $p_{n} = 0$, namely, if $x$ is a left
edge point, we let $h_{1} = gfg^{-n}$.  In all other cases we let
$h_{1} = g^{1-n}$.

The purpose of this is to reposition $x$ so that it is an interior
point of $C$.  In all cases, the point $h_{1}(x)$ has interior
position and sits in the interval $K_{1}$.

Choose two positions $p_{x}^{-}$ and $p_{x}^{+}$ from the set $P^{L}$
that ``bracket'' the point $h_{1}(x)$, namely:
$$(K_{1}, p_{x}^{-}) < h_{1}(x) < (K_{1}, p_{x}^{+})$$

We want to choose two positions $p_{y}^{-}$ and $p_{y}^{+}$ from
$P^{L}$ that bracket $y$ in a similar manner, and are also close, but
$y$ may be an edge point.  In that case, we can not bracket $y$, but
we can position ourselves close to it.

If $p_{m} = 1$, namely if $y$ is a right edge point, choose positions
such that
$$(K_{m}, 1) - (K_{m}, p_{y}^{-}) < \epsilon, p_{y}^{+} > p_{y}^{-}$$
Symmetrically, if $p_{m} = 1$, namely if $y$ is a left edge point,
choose positions such that
$$(K_{m}, p_{y}^{+}) - (K_{m}. 0) < \epsilon, p_{y}^{-} < p_{y}^{+}$$
In all other cases (namely, if $y$ is interior to $C$), choose
positions such that \mbox{$p_{y}^{-} < p_{m} < p_{y}^{+}$} and \mbox{$(K_{m},
p_{y}^{+}) - (K_{m}, p_{y}^{-}) < \epsilon$}.

The map $\pi$ is one-to-one and onto, and so there is some $r \in
\mathbb{N}$ such that $r = \pi((p_{x}^{-}, p_{x}^{+}, p_{y}^{-}, p_{y}^{+}))$.
Let $h_{2} = g^{m-r-1}fg^{r-1}$.

\begin{claim}
Let $h = h_{2}h_{1}$.  Then $|y - h(x)| < \epsilon$
\end{claim}

\emph{Proof:} For proof we first recall that $h_{1}(x)$ sits in
$K_{1}$, and has a position bracketed by $p_{1}$ and $p_{2}$.
All maps given are homeomorphisms that preserve orientation.  The
following diagram shows that $h$ satisfies Equation~\ref{eq:claim}.

\[ \begin{array}{lccccc}
  & (K_{1}, p_{x}^{-}) & < & h_{1}(x) & < & (K_{1}, p_{x}^{+}) \\
  g^{r-1}: & \downarrow &  & \downarrow & & \downarrow \\
  & (K_{r}, p_{x}^{-}) & < & g^{r-1}h_{1}(x) & < & (K_{r}, p_{x}^{+}) \\
  f: & \downarrow &  & \downarrow & & \downarrow \\
  & (K_{r+1}, p_{y}^{-}) & < & f(g^{r-1}(h_{1}(x))) & < & (K_{r+1}, p_{y}^{+}) \\
  g^{m-r-1}: & \downarrow & & \downarrow & & \downarrow \\
  & (K_{m}, p_{y}^{-}) & < & h(x) & < & (K_{m}, p_{y}^{+}) 
\end{array} \]

Whether $y$ is an edge or interior to $C$, we have ensured that $h(x)$
is squeezed within $\epsilon$ of it.

\subsection{Case 4: Level 2 Integer Type}  

\begin{figure}
\resizebox{\textwidth}{!}{\includegraphics{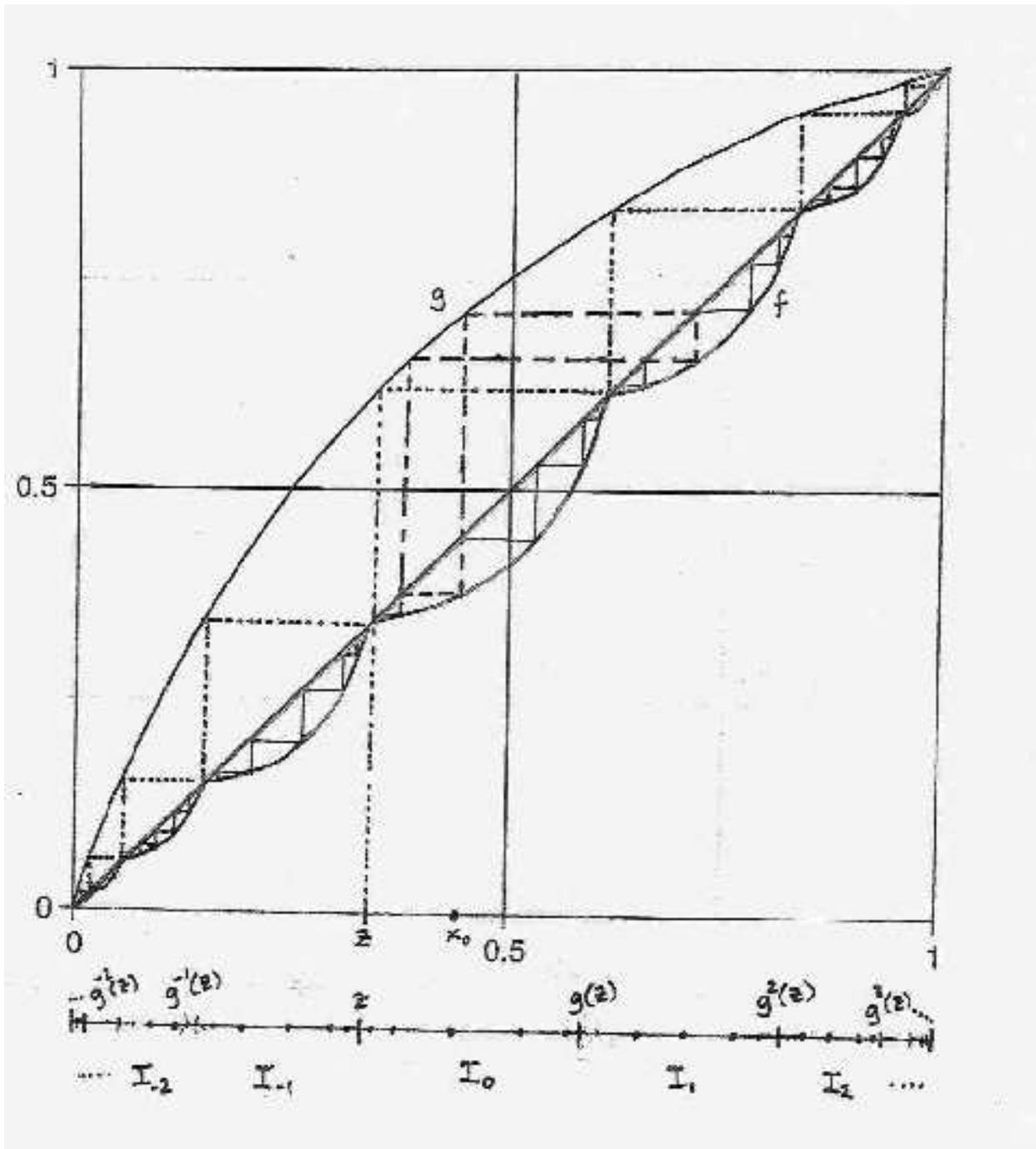}}
\caption{A level 2 integer type orbit.  The dotted line indicates the
level 1 integer type orbit of $z$.  The dashed line demonstrates the
commutivity of $f$ and $g$.  The full line follows the orbit of
$x_{0}$, which is recorded below the graph.}
\label{fig:lev2int}
\end{figure}

Let $g$ be a homeomorphism of the interval such that $g(x) > x$ for
all $x$ in $(0, 1)$.  Pick a point $z$ in the interior of the
interval.  We construct another map $f$ in such a way that the set
$\{g^{n}(z) : n \in \mathbb{Z}\}$ is preserved, so that $z$ is Case 2,
integer type.  The interaction between $f$ and $g$ off the orbit of
$z$ results in orbits that accumulate exactly along the closure of the
orbit of $z$, which we call \emph{level 2 integer type}. See Figure~\ref{fig:lev2int}.

We number the intervals which are the complement of \orb{z}: $I_{0}$
the interval to the right of $z$, $I_{n}$ the interval to the right of
$g^{n}(z)$. We write $f_{n}$ for the restriction of $f$ to $I_{n}$.

Let $f_{0}$ be a homeomorphism of $I_{0}$ with no interior fixed
points.  For all $n$, let
$$f_{n} = f|_{I_{n}} = g^{n}f_{0}g^{-n}(x), \, x \in I_{n}$$ 

Clearly $f$ and $g$ commute, $\orb{x_{0}} \cap I_{0} =
\{f^{n}(x_{0})\}$, and $\orb{x_{0}} \cap I_{n}$ is topologically
equivalent for all $n$.  Therefore $\orb{x_{0}}$ accumulates exactly
along $\overline{\orb{z}}$.  We call this \emph{level 2 integer type}.

We remark that by continuing to construct generators in this manner,
we could repeat this process and create points with orbits
accumulating exactly on the closure of the orbit of $x_{0}$,
\emph{level 3 integer type}, and orbits accumulating on that, etc.
See Section~\ref{sec:leveln} for a construction of a level $n$ integer
type orbit.

\subsection{Case 4: Level 1 Integer Type, Level 2 Dense}

\begin{figure}
\resizebox{\textwidth}{!}{\includegraphics{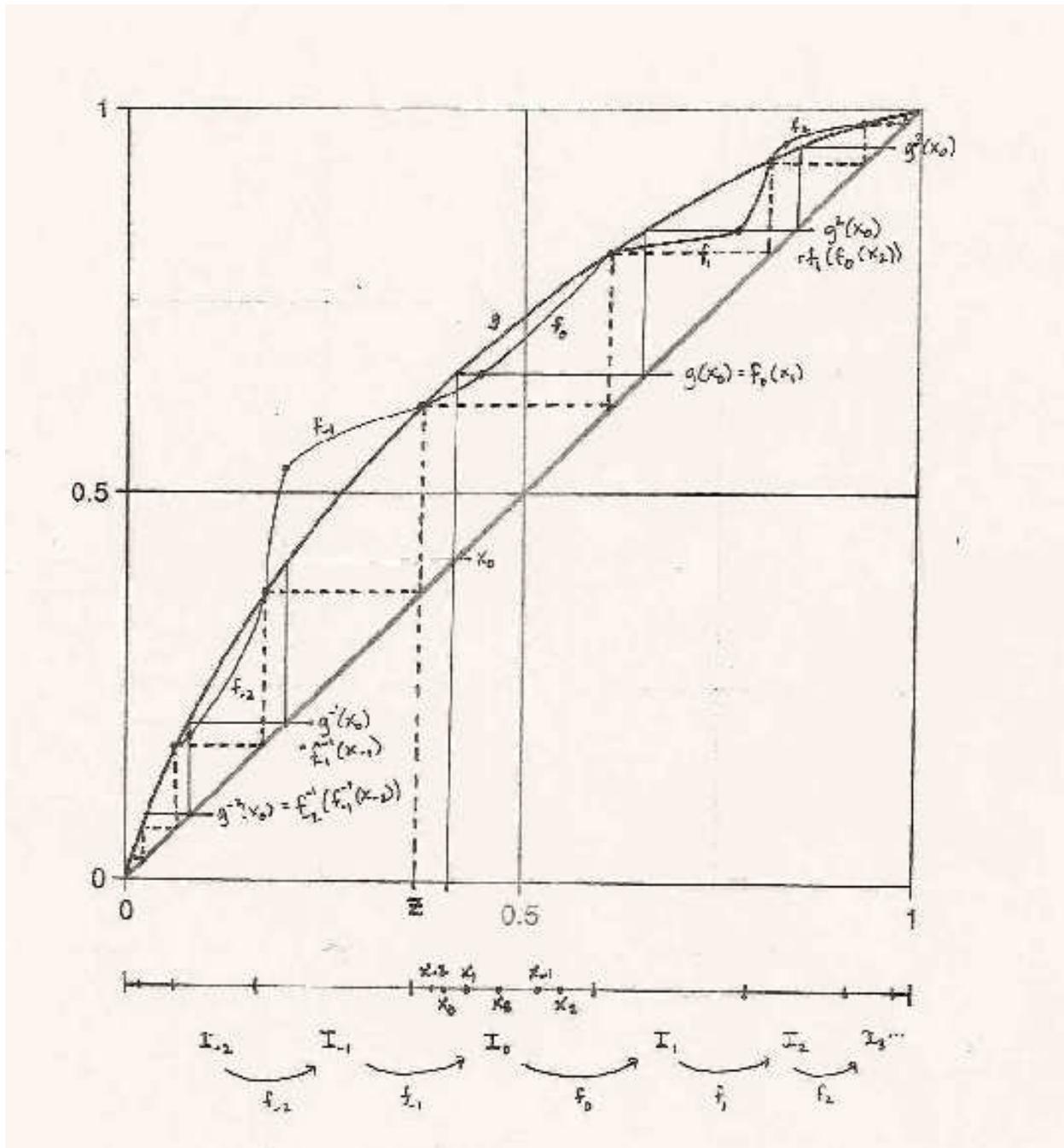}}
\caption{A dense level 2 orbit. The dashed line is the level 1 integer
type orbit of $z$.  The full line follows the action of $g$ on $x_{0}$, so that
$f$ can be defined accordingly, starting with $f_{0}$.}
\label{fig:lev2dense}
\end{figure}

As in the previous example, we start with a homeomorphism $g$ of the
interval with $g(x) > x$ in $(0,1)$, we choose some $z \in (0,1)$, and
let $I_{n}$ be the intervals complementary to $\{g^{n}(z) : z \in
\mathbb{Z}\}$.  Let $\{x_{n}\}$ be a countable set of points dense in
$I_{0}$. See Figure~\ref{fig:lev2dense} for a picture of the construction of $f$.

Let $f_{0}$ be a homeomorphism of $I_{0}$ to $I_{1}$ such that
$f_{0}(x_{1}) = g(x_{0})$.  This ensures $x_{1}$ is in the orbit of
$x_{0}$.  Next, let $f_{1}$ map $I_{1}$ to $I_{2}$ such that
$f_{1}(f_{0}(x_{2})) = g^{2}(x_{0})$.  This ensures that $x_{2}$ is in
the orbit of $x_{0}$.  Similarly, let
$f_{2}(f_{1}(f_{0}(x_{3})))=g^{3}(x_{0})$, and generally:

$$f_{n} \circ f_{n-1} \circ \cdots \circ f_{0}(x_{n+1}) = g^{n+1}(x_{0})$$

In defining each piece $f_{n}$, we map $I_{n} \rightarrow I_{n+1}$ and
specify the image of only one interior point.  The specified image
point depends only on $g$ and the pieces $f_{i}$ of $f$ we have
already defined.  We do a similar construction for $n < 0$.  

\orb{x_{0}} is dense in $I_{0}$ because it contains each $x_{n}$.
Since $g^{n}$ takes $I_{0}$ to $I_{n}$ homeomorphically, \orb{x} is
dense throughout the interval $\bar{I}$.  We have constructed a group
that has a level 1 integer type orbit and a level 2 dense orbit.

\subsection{Case 4: Level 1 Integer Type, Level 2 Cantor Set}

In order to use the clearest possible notation for this example, we
choose to think of the group as homeomorphisms of the real line,
rather than the interval.  The map $g$ will play the same role it has
in the previous examples: now we will think of it as translation by 1.
The level 1 integer type orbit will be $\orb{0} = \mathbb{Z}$.

Place a middle thirds Cantor set in each interval.  We denote by
$C_{n}$ the Cantor set in $[n, n+1]$.  The union of the Cantor sets is
invariant under $g$.  Let $\{x_{n}\}$ be the collection of left
endpoints of $C_{0}$.

Let $f_{0}$ be the homeomorphism of $[0,1]$ to $[1,2]$ with
$f_{0}(x_{1}) = g(x_{0}) = x_{0} + 1$ and $f_{0}(C_{0}) = C_{1}$ which
is given by Lemma~\ref{lemma:cantor1}.  Use the lemma again to find an
$f_{1}$ that will map $[1,2]$ to $[2,3]$, take $C_{1}$ to $C_{2}$, and
such that $f_{1}(f_{0}(x_{2})) = g^{2}(x_{0}) = x_{0} + 2$.
Similarly, let $f_{2}(f_{1}(f_{0}(x_{3})))=g^{3}(x_{0}) = x_{0} + 3$,
take $C_{2}$ to $C_{3}$, and generally, $f_{n}:I_{n} \rightarrow
I_{n+1}$ will take $C_{n}$ to $C_{n+1}$ and satisfy
$$f_{n} \circ f_{n-1} \circ \cdots \circ f_{0}(x_{n+1}) =
g^{n+1}(x_{0}) = x_{0} + (n+1)$$

The orbit of the origin is the integers; the Cantor set $\bigcup
C_{n}$ is invariant under the group. The point $x_{0}$ in the
Cantor set has an orbit which includes all the left endpoints of
$C_{0}$, making it dense in $C_{0}$, and by applying $g$, dense in
each $C_{n}$.  We say that $x_{0}$ has \emph{level 2 Cantor set type}.

\subsection{Level $n$}
\label{sec:leveln}

We have used the term \emph{level $n$} in the examples; we define it here.

\begin{defn}  
Orbits of points in $P(G)$ are level 0.  We say an orbit $\orb{z}$ is
\emph{level $n$} if it accumulates at the closure of a distinct
\emph{level $n-1$} orbit, and $n-1$ is the largest such integer.
\end{defn}

If $\orb{z}$ is in Case 1, 2 or 3 of the theorem, its accumulation at
$P(G)$ means it is level 1.  Higher level orbits arise in Case 4.  (We
have seen that it is possible for an orbit to be ``level 2 dense,''
and we consider this part of Case 4.)  We have the following:

\begin{theorem} 
There is a group on $n$ generators that admits a level $n$ orbit.
\end{theorem}

\begin{figure}
\resizebox{\textwidth}{!}{\includegraphics{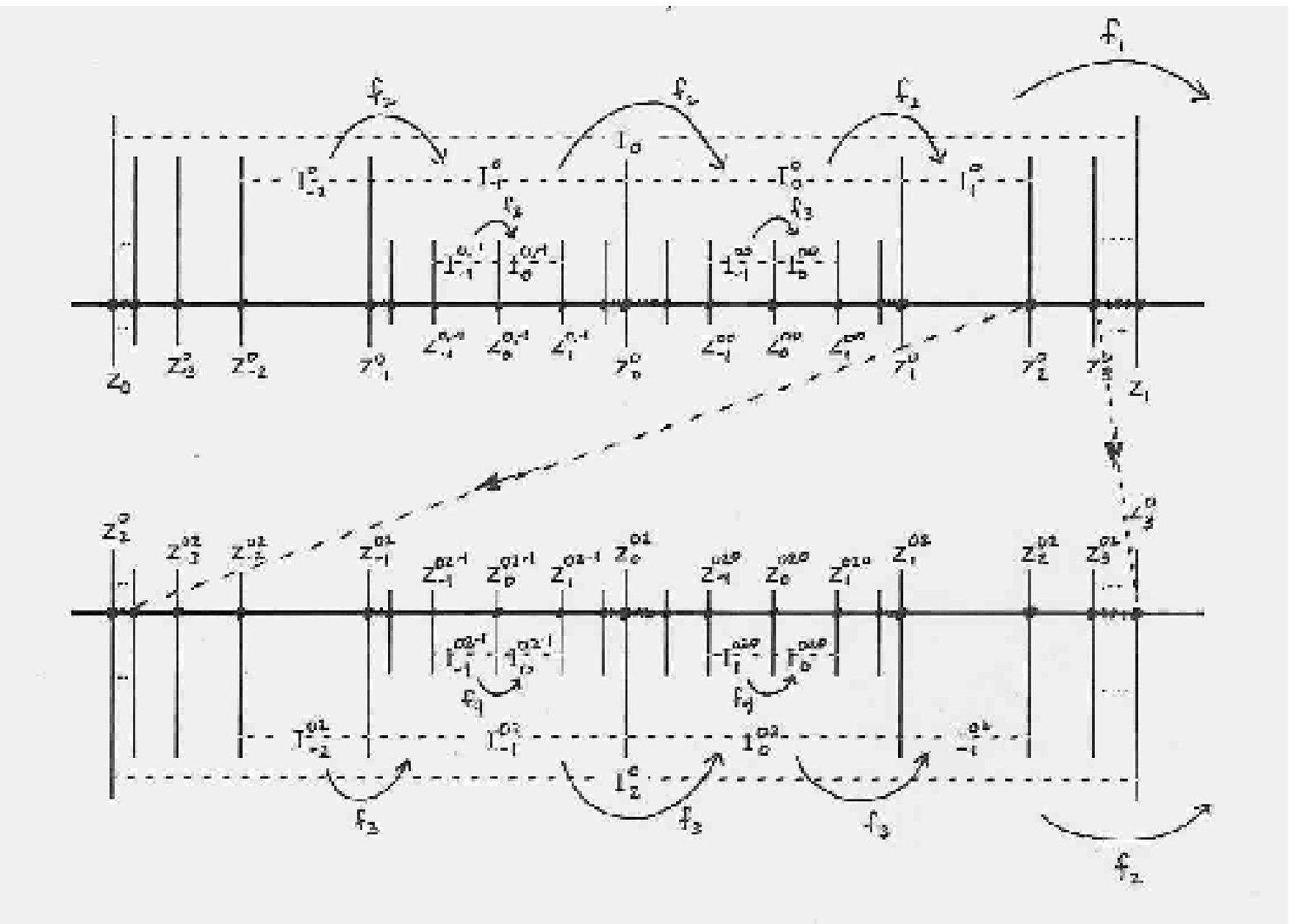}}
\caption{Notation for the Level n Integer Type Construction}
\label{fig:levnnotation}
\end{figure}

In the construction, we ``nest'' integer type orbits.  The notation
becomes cumbersome; however, with a picture in mind, we have tried to
make it as intuitive as possible.  See Figure~\ref{fig:levnnotation}.

\emph{Proof by construction:}

We base all our generators on a homeomorphism $f$ that takes an
interval $[a,b]$ to itself and satisfies $f(x) > x$ for all $x \in (a,
b)$.

Our first generator $f_{1}$ will be $f$ on $\bar{I} = [0,1]$.  We
choose a point $z_{0} \in (0,1)$ and we are careful that the remaining
generators preserve the orbit of $z_{0}$ under $f_{1}$, thereby
ensuring that $\orb{z_{0}}$ has level 1 integer type.  We write $z_{i}
= f_{1}^{i}(z_{0})$, and $\bar{I_{i}} = [z_{i}, z_{i+1}]$.  The
$I_{i}$ are the open subintervals complementary to the orbit of
$z_{0}$.  Note that $f_{1}:\bar{I_{i}} \rightarrow \overline{I_{i+1}}$.

Let our second generator $f_{2}$ be $f$ on $\bar{I_{0}}$, and for $x
\in \bar{I_{i}}$, $f_{2}(x) = f_{1}^{i}f_{2}f_{1}^{-i}(x)$.  Note:
$f_{2}$ commutes with $f_{1}$.

Next we choose a point $z^{0}_{0} \in I_{0}$.  This will be the point
with a level 2 integer type orbit.  We write $z^{0}_{i} =
f_{2}^{i}(z_{0}^{0})$, and $I_{i}^{0} = [z_{i}^{0}, z_{i+1}^{0}]$.
Notationally, the superscript $0$ tells us that the point is in
$I_{0}$, or that we have a subinterval of $I_{0}$.  We give notation
for the entire orbit of $z_{0}^{0}$ in the following way: Let
$z_{0}^{k} = f_{2}^{k}(z_{0}^{0})$, and $z_{i}^{k} =
f_{2}^{i}(z_{0}^{k})$.  The commutivity of $f_{1}$ and $f_{2}$ means
that we could just as well have defined $z_{i}^{0} =
f_{2}^{i}(z_{0}^{0})$ and $z_{i}^{k} = f_{1}^{k}(z_{i}^{0})$.
Similarly, we write $I_{i}^{k} = [z_{i}^{k}, z_{i+1}^{k}]$.  With this
notation we see that $f_{1}(z_{i}^{k}) = z_{i}^{k+1}$ and
$f_{2}(z_{i}^{k}) = z_{i+1}^{k}$, and more generally, $f_{2}:
I^{k}_{i} \rightarrow I^{k}_{i+1}$.  See
Figure~\ref{fig:levnnotation}.

As $f_{1}$ and $f_{2}$ commute, the orbit of $z_{0}^{0}$ is level 2
integer type, namely a countable collection of points in each $I_{i}$
accumulating at $\orb{z_{0}}$.  Occasionally, we will also call this a
\emph{second order nested} integer type orbit.

Let our third generator $f_{3}$ be $f$ on $I_{0}^{0}$.  We extend it
to all of $I_{0}$ by saying for $x \in I_{i}^{0}$, $f_{3}(x) =
f_{2}^{i}f_{3}f_{2}^{-i}(x)$, and from there to all of $I$ by saying
for $x \in I_{i}$, $f_{3}(x) = f_{1}^{i}f_{3}f_{1}^{-i}(x)$.  A little
work shows that $f_{3}$ commutes with both $f_{2}$ and $f_{1}$.

Choose $z_{0}^{00} \in I_{0}^{0}$.  This will be our point with a
level three orbit.  We write $z_{i}^{00} = f_{3}^{i}(z_{0}^{00})$ and
$I_{i}^{00} = [z^{00}_{i}, z^{00}_{i+1}]$.  Again, we extend this
notation to all of the orbit of $z_{0}^{00}$ in the following way: We
write $z_{i}^{k0} = f_{1}^{k}(z_{0}^{00})$ and $z_{i}^{kj} =
f_{2}^{j}(z_{i}^{k0})$.  Let $I_{i}^{kj} = [z_{i}^{kj}, z_{i+1}^{kj}]$.
Note that $f_{3}:I^{jk}_{i} \rightarrow I^{jk}_{i+1}$.

Notationally, at this point we can specify a point's position relative
to the level three orbit by noting which subinterval it is in: $x \in
I_{i}^{jk}$ is a point of the $i$th subinterval of the $k$th
subinterval of $I_{j}$, namely $x \in (z_{i}^{jk}, z_{i+1}^{jk})$.
The subscript refers to the most recent breakdown into subintervals,
and the superscripts place the point with respect to the previous
breakdowns (caused by the level one and two orbits).  Therefore a point with
two superscripts and one subscript is being specified with repsect to
the level three orbit.  

As we now wish to discuss the $n$th step, at which point the indices
begin to build up, we let $w(x,n)$ be a word of $n$ letters, each
letter an integer, specifying the location of the point $x$ with
respect to the level $n$ orbit.  By $w_{k}(x,n)$ we mean the $k$th
letter in the word $w(x,n)$.  When it is clear which word we are
referring to (as in the case where we are dealing only with a
particular point, or when we are using the word to specify all the
points in a particular subinterval), we drop reference to the point
$x$ and write merely $w(n)$ or $w_{k}(n)$.  We write $0(n)$ to
indicate the word of $n$ zeros.

Let our $n$th generator $f_{n}$ be $f$ on $I^{0(n-2)}_{0}$.  Extend
$f_{n}$ to $I^{0(n-2)}$: if $x \in I^{0(n-2)}_{i}$, $f_{n}(x) =
f_{n-1}^{i}f_{n}f_{n-1}^{-i}(x)$; and then extend to $I^{0(n-3)}$ by
shifting over and back by $f_{n-2}$, etc., all the way down the chain
of subintervals until $f_{n}$ is defined on all of $I$ (by shifting
over and back by $f_{1}$).  In this way $f_{n}$ preserves the orbits of
the previous level $k$ orbit points, $k < n$, and commutes with all
the $f_{k}$.  Choose $z_{0}^{0(n-1)} \in I^{0(n-2)}_{0}$.  We write
$I_{i}^{0(n-1)} = [z_{i}^{0(n-1)}, z_{i+1}^{0(n-1)}]$.

As before, for the complete level n orbit we write
$$z_{0}^{k 0(n-2)} = f_{1}^{k}(z_{0}^{0(n-1)})$$
$$z_{0}^{kj 0(n-3)} = f_{2}^{j}(z_{0}^{k 0(n-2)})$$
and generally
$$z_{0}^{w(k) 0(n-k-1)} = f_{k}^{w_{k}(k)}(z_{0}^{w(k-1) 0(n-k)})$$
$$z_{i}^{w(n-1)} = f_{n}^{1}(z_{0}^{w(n-1)})$$
and for the intervals, we write $I_{i}^{w(n-1)} = [z_{i}^{w(n-1)}, z_{i+1}^{w(n-1)}]$.
Note $f_{n}:I_{i}^{w(n-1)} \rightarrow I_{i+1}^{w(n-1)}$.

As $f_{n}$ acts like $f$ on $I_{0}^{0(n-2)}$, it gives
$z_{0}^{0(n-1)}$ an integer type picture in that interval, and the
orbit in that interval accumulates at $z_{0}^{0(n-2)}$ and
$z_{1}^{0(n-2)}$.  As all the generators commute, the orbit of
$z_{0}^{0(n-1)}$ is clearly an $n$th order nested integer type, and
has level $n$.

\pfbox

An interesting question comes from a variation on the converse of this
question; does the existence of a level $n$ orbit imply that the group
must have at least a certain number of generators?  We formulate a
conjecture along these lines in Section~\ref{sec:qsandconjs}.

\section{Semigroups}
\label{sec:semigroups}

We set aside groups for a moment and consider instead semigroups.  The
situation for semigroups is clearly different: because we are not
guaranteed that semigroup elements will have inverses, we can not
assume that all sets which are invariant under the action of the
semigroup will intersect any particular middle subinterval we choose,
as we did with groups in order to prove Lemma~\ref{lemma:zorn}.
Slightly different orbit structures can result from this, since we can
not expect that nonempty minimal sets will exist inside each
$\mathcal{O}^{+}(x)$.  We give an example to demonstrate one way in
which this can effect the orbit structure of points.

\subsection{The Idea}

We wish to construct an example of a semigroup such that certain
orbits do not fall into one of the four cases given for group orbits.
The general idea is the following:

Consider two maps, $f$ and $g$, taking $\bar{I}$ to itself.  Both maps
will have a single fixed point inside the interval, one nearer to 0
and one nearer to 1.  The fixed points will both be sources.  Then any
point to the outer side of either fixed point will only be able to
move farther towards the edge.  But points that sit between the two
fixed points will be moved in opposite directions by $f$ and $g$.  If
$f$ and $g$ are ``incommensurate'' in this area, they may create
orbits that are dense between the fixed points.


In order to give an example of such a semigroup, we will actually use
four maps: two that behave as described above, and two additional maps
which will allow points to shift back towards the middle section.  In
the middle of the interval, we will design $f$ and $g$ to look like
the example we gave for groups with a level 1 dense orbit.

\subsection{The Maps}

First we pick six points of $I$: $0< x_{1} < a_{0} < a_{1} < a_{2} <
a_{3} < x_{2} < 1$.  The $a_{n}$ should be evenly spaced: $a_{i+1} -
a_{i} = r$.  Let $I_{1} = [a_{0}, a_{1}], I_{2} = [a_{1}, a_{2}]$, and
$I_{3}=[a_{2}, a_{3}]$.  The maps will be defined piecewise over these
regions.

First, we give the two maps whose purpose is transportation. Note that
the desired effect of sending points back towards the middle is
achieved, as $h_{1}: I_{3} \rightarrow I_{2}$ and $h_{2}: I_{1}
\rightarrow I_{2}$ by strict translation.

\[ h_{1}(x) = \left\{ \begin{array}{ll}
                       \frac{a_{1}}{a_{2}}(x-a_{2}) + a_{1} & \mbox{if $x < a_{2}$} \\
		       x-r & \mbox{if $a_{2} \leq x \leq a_{3}$} \\
		       \frac{1-a_{2}}{1-a_{3}}(x-a_{3}) + a_{2} & \mbox{if $x > a_{3}$}
		       \end{array}
\right. \]

\[ h_{2}(x) = \left\{ \begin{array}{ll}
                      \frac{a_{1}-a_{0}}{a_{0}}(x-a_{0})+a_{1} & \mbox{if $x < a_{0}$} \\
		      x+r & \mbox{if $a_{0} \leq x \leq a_{1}$} \\
		      \frac{1-a_{2}}{1-a_{1}}(x-a_{1}) + a_{2} & \mbox{if $a_{1} < x$}
		      \end{array}
\right. \]

Next we give the two maps that represent the original concept.  The
important part is their action on $I_{2}$: $f: I_{2} \rightarrow
I_{3}$ acts like $x^{1/3}$ and $g:I_{2} \rightarrow I_{1}$ acts like
$x^{2}$.  The functions could be filled in however we liked in the
other sections, so long as $f$ fixed $x_{1}$ and $g$ fixed $x_{2}$.

\[ f(x) = \left\{ \begin{array}{ll}
        x_{1}(\frac{x}{x_{1}})^{2} & \mbox{if $x \leq x_{1}$} \\
	\frac{a_{2}-x_{1}}{a_{1}-x_{1}}(x-x_{1}) + x_{1} & \mbox{if $x_{1} < x < a_{1}$} \\
	r(\frac{x-a_{1}}{r})^{1/3} + a_{1} + r & \mbox{if $a_{1} \leq x \leq a_{2}$} \\
	\frac{1-a_{3}}{1-a_{2}}(x-1) + 1 & \mbox{if $a_{2} < x$}
	\end{array}
\right. \]

\[ g(x) = \left\{ \begin{array}{ll}
        \frac{a_{0}}{a_{1}}x & \mbox{if $x < a_{1}$} \\
	r(\frac{x-a_{1}}{r})^{2} + a_{0} & \mbox{if $a_{1} \leq x \leq a_{2}$} \\
	\frac{1-x_{2}-a_{1}}{1-a_{2}}(x-x_{2}) + x_{2} & \mbox{if $a_{2} < x < x_{2}$} \\
	(1-x_{2})(\frac{x}{1-x_{2}})^{1/3} + x_{2} & \mbox{if $x_{2} \leq x$}
	\end{array}
\right. \]

By combining these maps with the translations, we make the maps
$\hat{f} = h_{1}f$ and $\hat{g} = h_{2}g$, both of which take $I_{2}$
to itself.  These two maps commute; the details of this are
straighforward, but lengthy, and so they are omitted.

We let $\Gamma$ be the semigroup generated by $f, g, h_{1}$, and
$h_{2}$, and we denote the orbit of $x$ under the semigroup by
$\mathcal{O}^{+}(x)$ (the notation corresponds to that of the forward
orbit of $x$ under a full group).

It should seem intuitive now that $\hat{f}$ and $\hat{g}$ are
sufficiently incommensurate in the region $I_{2}$ that points in that
region can be moved densely throughout that region.

\begin{theorem} 
For $x \leq x_{1}$ (repsectively $x \geq x_{2}$), $\mathcal{O}^{+}(x)$
has integer type in $[0, x_{1}]$ (resp. $[x_{2}, 1]$), i.e. it is a
countable collection of points that accumulates exactly at $0$
(resp. $1$), and $\mathcal{O}^{+}(x) \cap (x_{1}, x_{2}) = \emptyset$.
For $a_{1} < x < a_{2}$, on the other hand, $\mathcal{O}^{+}(x) \cap
[0, x_{1}]$ and $\mathcal{O}^{+}(x) \cap [x_{2}, 1]$ are integer type,
but $\mathcal{O}^{+}(x) \cap [a_{1}, a_{2}]$ is dense.
\label{thm:semigroups}
\end{theorem}

The proof follows easily from the fact that $\ln{2}$ and $\ln{3}$ are
incommensurate.

\section{A Closer Look at Case 4}
\label{sec:case4}

We now return to full groups.  In order to further understand the
analysis of Case 4, we look at the intersection of the orbit and the
subinterval that sits between two consecutive points of a lower level
orbit.  We examine the subgroup of $G$ consisting of homeomorphisms
that preserve those subintervals $I_{n}$.  The intent is to restrict
ourselves to an area where we can apply Lemma~\ref{lemma:zorn} again,
and achieve the same conclusions as Theorem~\ref{thm:main}, with
$I_{n}$ playing the role of $I$, and the subgroup playing the role of
$G$.  Throughout this section, we assume that the chosen lower level
orbit $\orb{z}$ has integer or Cantor set type, and we use $I_{n}$ to
denote the intervals complementary to $\orb{z}$.  We will use $z$ and
$z_{0}$ interchangebly, depending on whether we mean to emphasize that
it is the endpoint of the interval $I_{0}$.  We have a few lemmas.

\subsection{Some Useful Lemmas}

Let $H_{n} = \{g \in G : g(I_{n}) = I_{n}\}$.  This is clearly a
subgroup of $G$.  Let $\mathcal{O}_{H}(x) = \{ h(x) : h \in H\}$.

\begin{lemma} 
If $y \in I_{n}$, then $\orb{y} \cap I_{n} = \mathcal{O}_{H_{n}}(y)$.
\end{lemma}

\emph{Proof:} Let $y_{0} \in \orb{y} \cap I_{n}$.  So $y_{0} \in
I_{n}$ and for some $g_{0} \in G$, $y_{0} = g_{0}(y)$.  Group elements
must take complementary intervals to other complementary intervals,
so this means $g_{0} \in H_{n}$, and $y_{0} \in \mathcal{O}_{H_{n}}(y)$.

Let $y_{0} \in \mathcal{O}_{H_{n}}(y)$.  Then $y_{0} = g_{0}(y)$ for
some $g_{0} \in H_{n}$.  Clearly $y_{0} \in \orb{y} \cap I_{n}$.

\pfbox

Therefore, in order to get a picture of the group orbit of a point in
a given subinterval, it is enough to look at the orbit of the point under
the subgroup of homeomorphisms sending the interval to itself.

\begin{lemma}
\label{lemma:map_intervals}
If $\orb{z}$ has integer type, for all pairs of integers $i, j$, there
exists some $g \in G$ such that $g(I_{i}) = I_{j}$.
\end{lemma}

\emph{Proof:} Every group element $g$ must take intervals to
other intervals, As the maps are all
continuous, $g(I_{i})$ is totally determined by where $g$ takes an
endpoint $z_{i}$.  As each endpoint is part of the orbit of $z
= z_{0}$, there is some $g_{1} \in G$ with $g_{1}(z_{0}) = z_{i}$ and
some $g_{2} \in G$ with $g_{2}(z_{0}) = z_{j}$.  Therefore
$g_{2}g_{1}^{-1}(z_{i}) = z_{j}$ and so $g_{2}g_{1}^{-1}(I_{i}) =
I_{j}$.  

\pfbox

\begin{lemma} 
If $\orb{z}$ has integer type, then for all $i$ and $j, H_{i} =
H_{j}.$
\end{lemma}

\emph{Proof:} Suppose $h \in G$ takes $I_{n}$ to $I_{m}$.  Then by the
fact that $h$ is continuous and the intervals are adjacent, it is
clear that for all $k \in \mathbb{Z}$, $h$ must take $I_{n+k}$ to
$I_{m+k}$.  Chose $k$ such that $i + k = j$. Then $h \in H_{j}$ as
well.  The argument is symmetric, so $H_{i} = H_{j}$.

\pfbox

In this case ($\orb{z}$ of integer type) we shall refer to the
subgroup of elements that take each $I_{i}$ to itself as
$H^{1}$, and we have:

\begin{lemma} 
If there exists a level 2 orbit, $H^{1}$ is nontrivial.
\end{lemma}

\emph{Proof:} Suppose $y$ is a level 2 orbit point and $y \in
I_{0}$.  $\orb{y}$ accumulates along $\orb{z}$, so
there are infinitely many points of $\orb{y}$ in $I_{0}$.  Pick some
$y_{k} \neq y$; then $y_{k} = g(y)$ for some group element $g$.  That
group element is clearly in $H^{1}$ as it takes a point of $I_{0}$ to
another point of $I_{0}$, and it is nontrivial as $y_{k} \neq y$.

\pfbox

In the case that $\orb{z}$ has Cantor set type, it is not immediately
clear that each subinterval is a homeomorphic copy of $I_{0}$.
Because the orbit of $z$ is merely dense in the Cantor set, but not
necessarily equal to the Cantor set, we can not always find a map
taking $z$ to another particular endpoint.  However, the group orbit
will look the same in those subintervals for which there is such a
map.

\subsection{Parallel Orbits and Subdivisions}

\begin{defn}
If $I_{n}$ are the intervals making up the complement of $\orb{z}$,
and if for all $n, \orb{x} \cap I_{n}$ consists of one point, we say
that $\orb{x}$ is \emph{parallel} to \orb{z}.
\end{defn}

\begin{lemma}
If $\orb{z}$ has integer type, then any parallel orbit is also integer
type, and has the same level.  
\end{lemma}

If $\orb{z}$ has Cantor set type, what does a parallel orbit look
like?  The orbit consists of a single point in each open interval that
is part of the complement of the Cantor set.  We know the orbit
accumulates along the endpoints of the Cantor set, since it is part of
the range of the Cantor set orbit. In this case it does so in a
one-sided manner, by which we mean that the accumulation set at any
particular endpoint of an interval $I_{n}$ will consist of points
which are not contained in $I_{n}$, but sit exclusively to the other
side of the endpoint.  Therefore the orbit of $x$ is not integer type,
but neither is it inside the Cantor set; the closure looks like the a
Cantor set and the union of a countable collection of isolated points.
Unlike an orbit parallel to an integer type orbit, this orbit will
have a higher level, since it will, by necessity, accumulate at the
Cantor set to which it is parallel.

\emph{Remark:} The above description demonstrates that the concept
\emph{parallel} is not symmetric: $\orb{x}$ may be parallel to
$\orb{y}$ without $\orb{y}$ being parallel to $\orb{x}$.

For example, let $g$ be a homeomorphism with $g(x) > x$ for $x$ in
$(0,1)$.  Let $z$ be a point with integer type orbit, and $I_{n}$ the
complementary intervals.  Choose a point $y \in I_{0}$.  We let $f(x)
= g(x)$ for $x \in I_{n}, x \leq g^{n}(y)$, and for $x \in I_{n}, x >
g^{n}(y)$, let $f$ be any homeomorphism.  Now all $x \in I_{n}, x \leq
g^{n}(y)$ will have an integer type orbit which intersects each
$I_{n}$ in exactly one point.  Thus we have a whole interval of orbits
parallel to that of $z$.  Points of $I_{n}$ with $x > g^{n}(y)$ will
have more complicated orbits, depending on the interaction of $f$ and
$g$ in that areas.


\begin{lemma}
If $\orb{z}$ is integer type, and if $\orb{y} \cap I_{0}$ contains
exactly one point, then for all $n$, $\orb{y} \cap I_{n}$ contains
exactly one point, i.e., $\orb{y}$ is parallel to $\orb{z}$.
\end{lemma}

\emph{Proof:} Fix $n$.  From Lemma~\ref{lemma:map_intervals} we have
some $g \in G$ taking $I_{0}$ to $I_{n}$.  It is a homeomorphism, so
$\orb{y} \cap I_{n}$ contains exactly one point as well.

\pfbox

It is useful to note that if $\orb{y}$ accumulates at $\orb{z}$, and
$\orb{z}$ has integer type, then any orbits parallel to $\orb{z}$ must
sit exclusively to one side or the other of $\orb{y}$.  In other
words, if $y \in I_{0}$ and $P = \{x \in I: \orb{x}$ is parallel to
$\orb{z}\}$, then the set $P \cap I_{0}$ is contained completely in
$(z_{0}, y)$ or $(y, z_{1})$.  This is because any element $h \in
H^{1}$ will fix the (single) point of $\orb{x} \cap I_{0}$ and send
each half of the interval to itself.  Therefore, a single parallel
orbit relegates the orbit of $y$ to one side or the other of $I_{0}$,
and if follows, of each $I_{n}$.  In order for $\orb{y}$ to accumulate
at $\orb{z}$, it must not be restricted to an interior subinterval.
(If $\orb{z}$ were of Cantor set type, the complementary intervals
would not be adjacent, and a parallel orbit would not restrict other
points from accumulating along $\orb{z}$.  In fact, the orbits of all
points complementary to the Cantor set \emph{will} accumulate along
the Cantor set: see the proof of Theorem~\ref{thm:main}) This means
all parallel integer type orbits will sit to the same side of $y$, and
to the same side of all the points of $\orb{y} \cap I_{0}$.  In this
vein, when considering the behavior of $\orb{y}$ in an interval
$I_{0}$, it is reasonable to restrict to a subinterval of $I_{0}$ in
which there are no parallel orbits to $z$: for instance, $(x, z_{1})$
will be such an interval if $x$ is the rightmost parallel orbit point
in $I_{0}$.  We call this subinterval $J_{0}$, and its homeomorphic
images in each $I_{n}$ we call $J_{n}$.  Note that elements of $H^{1}$
will also take each $J_{n}$ to itself.

\begin{lemma}
$\orb{y} \cap I_{0} = \orbh{y}$, and therefore the group orbit of
$y$ intersected with $I_{0}$ can not be topologically different from the
that of $y$ intersected with any other $I_{i}$.
\end{lemma}

\emph{Proof:} Let $p \in \orb{y} \cap I_{0}$.  The equality is given
by Lemma~\ref{lemma:map_intervals}.  As the orbit of $y$ intersected
with any given $I_{i}$ contains the homeomorphic image of the orbit of
$y$ intersected with $I_{0}$, and vice versa, the two pictures are
topologically equivalent.  

\pfbox

We know that it is impossible for $y$ to have a Cantor set orbit in
one $I_{i}$ and an integer type orbit in another, or to be dense in
one and contained in a Cantor set in another.  However, the orbit of
$y$ may be different in seperate subintervals $J_{i}$ of $I_{n}$,
bounded by points of various parallel orbits, as elements of $H^{1}$
may behave in different ways on different subintervals of a given
$I_{n}$.

For both the integer type and the Cantor set type, we can restrict
ourselves to looking at a single $I_{n}$, bounded by consecutive
points of some lower level orbit, and develop the picture there.  (In
the integer case, we need only look at one subinterval, because the
orbit will be the same in all subintervals.)  If the lower level orbit
has orbits which are parallel to it, we may want to break the $I_{n}$
down into subintervals bounded by points of parallel orbits.  We
expect that in each subinterval, we have a repetition of the main
picture: dense orbits, a Cantor set, integer type, and possibly higher
levels.  However, in order to conclude that our restriction again
falls into one of the original four cases, we need to be able to apply
Lemma~\ref{lemma:zorn} to the subgroup of $G$ associated to that
interval.  Unfortunately, we do not know if either of the conditions
are satisfied: subgroups of finitely generated groups are not
necessarily themselves finitely generated, and we do not have enough
information about the maps in the subgroups to conclude that one of
them is free of accumulating fixed points.  In order to conclude
anything additional, we therefore require additional hypotheses.

\section{Analytic Groups}
\label{sec:analytic}

\begin{theorem}
\label{thm:analytic}
If $G$ is a finitely generated group of homeomorphisms of the interval
$\bar{I}$ which are continuous on $\bar{I}$ and analytic on $I$, and
if there are no fixed points of the group in $I$, then there can be at
most one integer type orbit with level greater than 1.  In addition,
if $\orb{x}$ is a level $k$ integer type orbit, and it intersects some
$I_{n}$ (an interval complementary to the level $k-1$ orbit), all
other orbits intersecting $I_{n}$ will be parallel to $\orb{x}$, i.e.,
they will all be level $k$ integer type orbits.  Therefore there will
be no orbits of level greater than $k$ in $I_{n}$ and its homeomorphic
copies.
\end{theorem}

\emph{Proof:}

Since $G$ is finitely generated and there are no fixed points of the
group inside $I$, we can apply Lemma~\ref{lemma:zorn} to find a
$G_{I}$-minimal set inside $\overline{\orb{x}}$.  As in
Theorem~\ref{thm:main}, we can apply minimality arguments to conclude
that $\orb{x}$ must fall into one of the four cases.  Let us assume it
is a level $k$ integer type orbit, where $k$ is greater than 1.  Let
$\orb{z}$ be the level $k-1$ orbit on which $\orb{x}$ accumulates.
Since $\orb{x}$ has integer type, $\orb{z}$ must have integer or
Cantor set type. In either case, let $I_{n}$ be the collection of
complementary intervals, such that $x \in I_{0}$.  Since $x$ has level
$k$ integer type, we can further divide $I_{0}$ into subintervals
complementary to the orbit of $x$, which we call $I_{0n}$.

Consider the subgroup $H_{0} = \{g \in G: g(I_{0}) = I_{0}\}$.
$H_{0}$ is clearly nontrivial as $\orb{x}$ intersects $I_{0}$ in an
infinite number of points.  But what about $H_{00} = \{h \in H_{0} :
h(I_{00}) = I_{00}\}$?  We note that $H_{00} = H_{0n}$ for all $n \in
\mathbb{Z}$, since the $I_{0n}$ are adjacent.  But this means that $h$
has a sequence of fixed points accumulating at the endpoints of
$I_{0}$, namely, the endpoints of the $I_{0n}$.  Since $I_{0}$ is
interior to $I$, and $h$ is analytic in $I$, this means $h$ must be
the identity map.  Therefore every point $y \in I_{00}$ has an orbit
which intersects each $I_{0n}$ in exactly one point, giving it an
orbit parallel to that of $x$, level $k$ integer type.  Thus $k$ is
the highest order level for orbits of points in $I_{0}$.

\pfbox

\begin{theorem}
If $G$ is a finitely generated abelian group of homeomorphisms of
$\bar{I}$ which are continuous on $\bar{I}$ and analytic on $I$, there
can be no orbits of Cantor set type.  Therefore the highest possible
non-dense orbit is level 2 integer type.
\end{theorem}

\emph{Proof:}

Suppose $\orb{x}$ were Cantor set type.  Let $I_{n}$ be the collection
of intervals complementary to the Cantor set, and say $x_{0}^{1} \in
\orb{x}$ is the left endpoint of the interval $I_{0}$.  Consider
$H_{n} = \{ g \in G: g(I_{n}) = I_{n}\}$.  The orbit of $x_{0}^{1}$ is
dense in the Cantor set, so if $x_{n}^{1}$ is the left endpoint of
$I_{n}$, there must be a sequence of maps $g_{k}$ such that $x_{k}^{1}
= g_{k}(x_{0}^{1}) \rightarrow x_{n}^{1}$, and $I_{0}$ is mapped to
the corresponding intervals $I_{k}$.  Consider $H_{k}$.  Since
$g_{k}:I_{0} \rightarrow I_{k}$, then if $h \in H_{0}$, $ghg^{-1} \in
H_{k}$.  Because $G$ is abelian, this means $h \in H_{k}$.  So $h$ has
each $x_{k}^{1}$ as a fixed point, and the $x_{k}^{1}$ accumulate at
$x_{n}^{1}$.  However, $h$ is analytic in $I$, and $x_{n}^{1}$ is interior
to $I$.  This is a contradiction.  So $\orb{x}$ must not be Cantor set
type.  Therefore, if $\orb{x}$ is not dense, it must be either level 1
integer type, or level 2 integer type (accumulating along a level 1
integer type).  It can not be level 3, by Theorem~\ref{thm:analytic}.

\pfbox

\section{Questions and Conjectures}
\label{sec:qsandconjs}

We would like to be able to make a definitive statement about the
overall structure of orbits of points in Case 4.  It seems that the
natural thing to expect is that orbits in Cases 1, 2, or 3 ``nest'' at
increasing levels, such that the orbit of a level $n$ point $x$ can be
understood by first looking at a single interval, broken down into
subintervals bounded by consecutive parallel orbits, such that in each
section the orbit will be dense, integer type, or contained in a
Cantor set.  That picture will then be repeated a countable number of
times at each level.  However, in order to conclude this, we must
better understand the structure of the subgroups that preserve
subintervals at a given level.  It seems likely that some amount of
additional smoothness of the generators will enable us to apply the
crucial lemma to the subgroups.  We expect that we need something
stronger than continuity, but that it is not necessary to require as
much as analyticity:

\begin{conjecture}
If $G$ is a group generated by $C^{2}$ homeomorphisms, orbits of
points in Case 4 will consist of a countable number of copies of an
interval $I$, the union of a number of subintervals $J_{n}$, in each
of which the orbit is dense, has integer type, or is contained in a
Cantor set.
\end{conjecture}

Without knowing exactly how much smoothness is required, it would be
interesting to see an example of a group of homeomorphisms with a
level 1 integer type orbit where the subgroup $H^{1}$ that preserves
the subintervals does not satisfy either condition in
Lemma~\ref{lemma:zorn}.  One imagines that the generating
homeomorphisms of such a group would be quite complicated, either
individually or in their interaction.

It seems likely that more statements could be made about the
level orbits that are possible under a particular group:

\begin{conjecture}
If $G$ is a group generated by $n$ homeomorphisms, the highest
possible level of an orbit is $n+1$.
\end{conjecture}

There are also natural questions about the effect of the group
structure on the possible orbits.  We saw that the addition of abelian
to the analytic case greatly restricted the possibilities.  In
particular, if the group is abelian, it may simplify the complications
created by an orbit accumulating along a Cantor set.  What effect
might abelian have on a group of less smooth homeomorphisms?  What
about other group structure restrictions?  This issue is partially
addressed in \cite{FF}.

Finally, in the case of semigroups, although the structures can
clearly be different, there is probably a similar description for
orbits which takes into account any sections of the circle where all
maps move points in the same direction.


\begin{thebibliography}{00}

\bibitem{FF}
Benson Farb and John Franks, Groups of Homeomorphisms of One-Manifolds III: Nilpotent Subgroups. \emph{Ergodic Theory and Dynam. Syst.}, 2001

\bibitem{ghys}
Etienne Ghys, Groups Acting on the Circle, IMCA, Lima, 1999

\bibitem{katok}
Anatole Katok and Boris Hasselblatt, \emph{Introduciton to the Modern Theory of Dynamical Systems}, Cambridge University Press, New York, 1995

\bibitem{pugh}
Charles Pugh, \emph{Introduction to Analysis}, Springer, New York, 2000

\bibitem{robin}
Clark Robinson, \emph{Dynamical Systems: Stability, Symbolic Dynamics, and Chaos}, CRC Press, Boca Raton, 1998

\end{thebibliography}

\end{document}